\documentclass[final,1p,times,twocolumn]{elsarticle}

\usepackage{natbib}

 \usepackage{graphicx}
\usepackage{epstopdf}

\usepackage{subfigure}

\usepackage{latexsym,bm}

\usepackage{amssymb}
\usepackage{amsmath}
\usepackage{amsfonts}

\newtheorem{thm}{Theorem}
\newtheorem{lem}{Lemma}

\newtheorem{rem}{Remark}

\newcommand{\norm}[1]{\left\Vert#1\right\Vert}

\newcommand{\pref}[1]{(\ref{#1})}

\allowdisplaybreaks

\begin{document}

\begin{frontmatter}

\title{The Convergence of Least-Squares Progressive Iterative Approximation with Singular Iterative Matrix}



\author[linadd]{Hongwei Lin\corref{cor1}}
\author[linadd]{Qi Cao}
\author[linadd]{Xiaoting Zhang}
\cortext[cor1]{Corresponding author (hwlin@zju.edu.cn). }
\address[linadd]{School of Mathematical Science, State Key Lab. of CAD\&CG, Zhejiang University, Hangzhou, 310027, China}

\date{}

\begin{abstract}
 Developed in [Deng and Lin, 2014],
    Least-Squares Progressive Iterative Approximation (LSPIA) is an efficient iterative method for solving B-spline curve and surface least-squares fitting systems.
 In [Deng and Lin 2014],
    it was shown that LSPIA is convergent when the iterative matrix is nonsingular.
 In this paper, we will show that LSPIA is still convergent even the
    iterative matrix is singular. 
\end{abstract}


\end{frontmatter}



\section{Introduction}

 Least-squares fitting is a commonly employed approach in
    engineering applications and scientific research,
    including geometric modeling.
 With the advent of \emph{big data} era, least-squares fitting systems with
    singular coefficient matrices often appear,
    when the number of the fitted data points is very large,
    or there are ``holes'' in the fitted data points.
 LSPIA~\cite{deng2014progressive} is an efficient iterative method for
    least-squares B-spline curve and surface fitting~\cite{brandt2015optimal}.
 In Ref.~\cite{deng2014progressive}, it was shown that LSPIA is convergent
    when the iterative matrix is \emph{nonsingular}.
 In this paper, we will show that,
    when the iterative matrix is \emph{singular},
    LSPIA is still convergent.
 This property of LSPIA will promote its applications in large scale data fitting.

 The motivation of this paper comes from our research practices,
    where some singular least-squares fitting systems emerge.
 For examples, in generating trivariate B-spline solids by fitting
    tetrahedral meshes~\cite{lin2015constructing},
    and in fitting images with holes by T-spline surfaces~\cite{lin2013efficient},
    coefficient matrices of least-squares fitting systems are singular.
 There, LSPIA was employed to solve the least-squares fitting systems,
    and converged to stable solutions.
 However, in Ref.~\cite{lin2015constructing,lin2013efficient},
    convergence of LSPIA for solving singular linear systems was not proved.

 The progressive-iterative approximation (PIA) method was first developed
    in~\citep{lin2004constructing,lin2005totally},
    which endows iterative methods with geometric meanings,
    so it is suitable to handle geometric problems appearing in the field of geometric design.
 It was proved that the PIA method is convergent for B-spline
    fitting~\citep{lin2011extended,deng2014progressive},
    NURBS fitting~\citep{shi06iterative}, T-spline fitting~\citep{lin2013efficient}, subdivision surface fitting~\citep{cheng2009loop, fan2008subdivision,
    chen2008progressive}, as well as curve and surface fitting with totally positive basis~\citep{lin2005totally}.
 The iterative format of geometric interpolation
    (GI)~\cite{maekawa2007interpolation} is similar as that of
    PIA.
 While PIA depends on the parametric distance,
    the iterations of GI rely on the geometric distance.
 Moreover, the PIA and GI methods have been employed in some applications,
    such as reverse engineering~\citep{kineri2012b,yoshihara2012topologically}, curve design~\citep{okaniwa2012uniform}, surface-surface intersection~\citep{lin2014affine},
    and trivariate B-spline solid generation~\citep{lin2015constructing}, etc.

 The structure of this paper is as follows.
 In Section~\ref{sec:it_fmt}, we show the convergence of LSPIA with singular
    iterative matrix.
 In Section~\ref{sec:example}, an example is illustrated.
 Finally, Section~\ref{sec:conclusion} concludes the paper.

\section{The iterative format and its convergence analysis}
\label{sec:it_fmt}

 To integrate the LSPIA iterative formats for B-spline curves,
    B-spline patches,
    trivariate B-spline solids, and T-splines,
    their representations are rewritten as the following form,
    \begin{equation} \label{eq:unified_form}
      \bm{P}(\bm{t}) = \sum_{i=0}^n \bm{P}_i B_i(\bm{t}).
    \end{equation}
 Specifically, T-spline patches~\cite{sederberg2004t} and trivariate T-spline
    solids~\cite{zhang2012solid} can be
    naturally written as the form~\pref{eq:unified_form}.
 Moreover,
 \begin{itemize}
   \item If $\bm{P}(\bm{t})$~\pref{eq:unified_form} is a B-spline curve,
            then, $\bm{t}$ is a scalar $u$, and $B_i(\bm{t}) = N_i(u)$,
            where $N_i(u)$ is a B-spline basis function.
   \item If $\bm{P}(\bm{t})$~\pref{eq:unified_form} is a B-spline patch
            with $ n_u \times n_v$ control points,
            then, $\bm{t} = (u,v)$, and $B_i(\bm{t}) = N_i(u) N_i(v)$,
            where $N_i(u)$ and $N_i(v)$ are B-spline basis functions.
         In the control net of the B-spline patch,
            the original index of $N_i(u)$  is $[\frac{i}{n_u}]$,
            and the original index of $N_i(v)$ is $(i\ \text{mod}\ n_u)$,
            where $[\frac{i}{n_u}]$ represents the maximum integer not exceeding $\frac{i}{n_u}$,
            and $(i\ \text{mod}\ n_u)$ is the module of $i$ by $n_u$.
   \item If $\bm{P}(\bm{t})$ is a trivariate B-spline solid with
            $n_u \times n_v \times n_w$ control points,
            then $\bm{t} = (u,v,w)$, and $B_i(\bm{t}) = N_i(u) N_i(v) N_i(w)$.
         In the control net of the trivariate B-spline solid,
            the original index of $N_i(w)$ is $[\frac{i}{n_u n_v}]$,
            the original index of $N_i(u)$ is $[\frac{(i\ \text{mod}\ n_u n_v)}{n_u}]$,
            and the original index of $N_i(v)$ is
            $((i\ \text{mod}\ n_u n_v)\ \text{mod}\ n_u)$.
 \end{itemize}

\begin{figure}[!htb]
  \centering
    \label{subfig:valley}
    \includegraphics[width=0.6\textwidth]{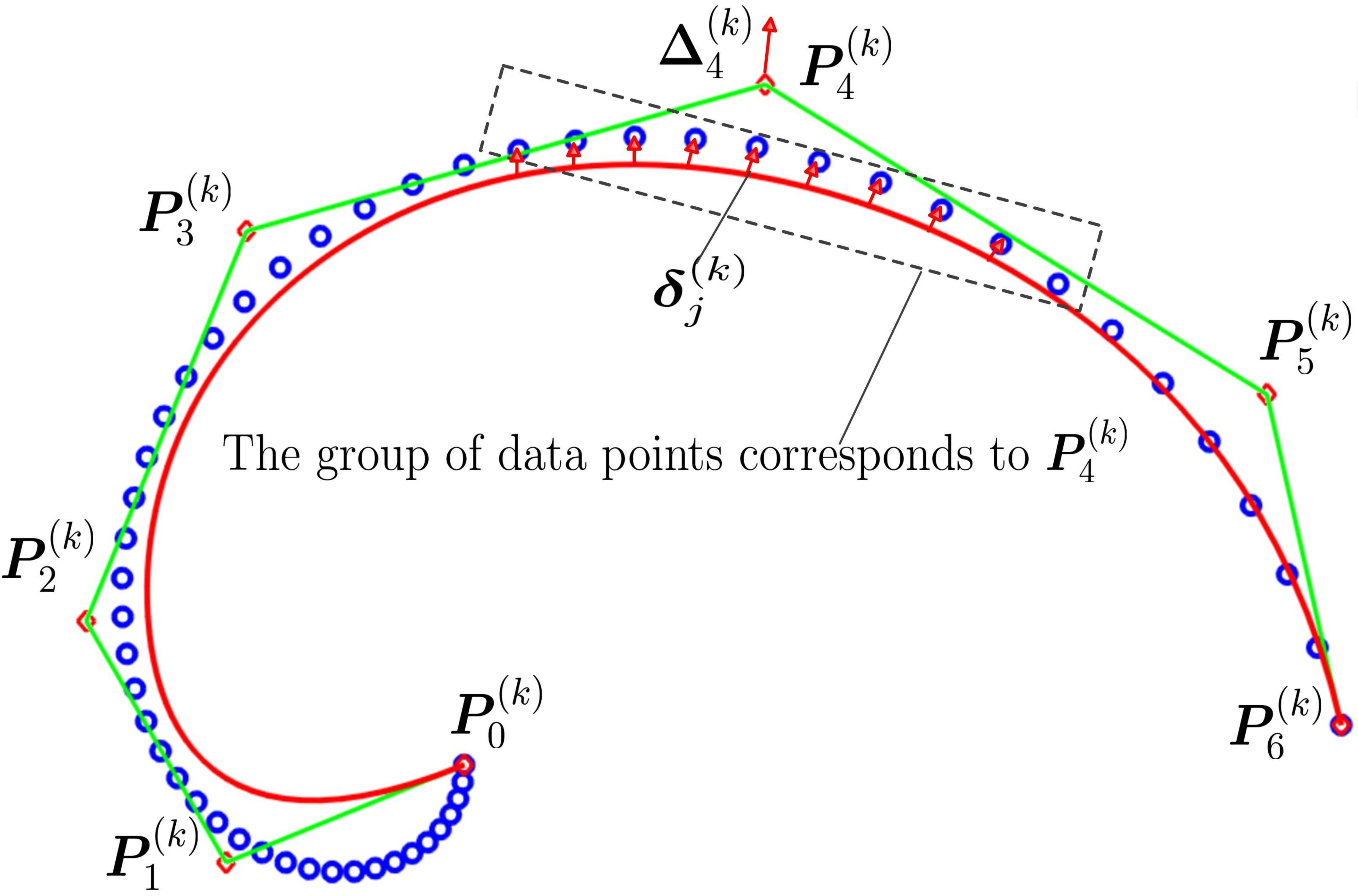}
    \caption{One iteration step of LSPIA includes two procedures,
        vector distribution and vector gathering.
    In the vector distribution procedure,
        all of DVDs $\bm{\delta}_j^{(k)}$ corresponding to a group of data points are distributed to the control point the data point group corresponds to.
    In the vector gathering procedure,
        all of DVDs distributed to a control point are weighted averaged to generate the DVC $\bm{\Delta}_i^{(k)}$.
    Here, blue circles are the data points,
        and the red curve is the $k^{th}$ curve $\bm{P}^{(k)}(u)$.}
   \label{fig:lspia_sketch}
\end{figure}

 Suppose we are given a data point set
    \begin{equation} \label{eq:data_points}
        \{\bm{Q}_j = (x_i,y_i,z_i), j = 0,1,\cdots, m\},
    \end{equation}
    each of which is assigned a parameter $\bm{t}_j, j = 0,1,\cdots,m$.
 Let the initial form be,
    \begin{equation} \label{eq:initial_form}
      \bm{P}^{(0)}(\bm{t}) = \sum_{i=0}^n \bm{P}_i^{(0)} B_i(\bm{t}),
      \ n \leq m.
    \end{equation}
 It should be noted that,
    though the initial control points $\bm{P}_i^{(0)}$ are usually chosen from the given data points,
    the initial control points are unrelated to the convergence of LSPIA.
 To perform LSPIA iterations,
    data points are classified into groups.
 All of data points with parameters $\bm{t}_j$ satisfying
    $B_i(\bm{t}_j) \neq 0$ are classified into the $i^{th}$ group,
    corresponding to the $i^{th}$ control point~\pref{eq:initial_form}.

 After the $k^{th}$ iteration of the LSPIA,
    the $k^{th}$ form $\bm{P}^{(k)}(\bm{t})$ is generated,
    \begin{equation*}
      \bm{P}^{(k)}(\bm{t}) = \sum_{i=0}^n \bm{P}_i^{(k)} B_i(\bm{t}).
    \end{equation*}
 To produce the $(k+1)^{st}$ form $\bm{P}^{(k+1)}(\bm{t})$,
    we first calculate the \emph{difference vectors for data points} (DVD) (Fig.~\ref{fig:lspia_sketch}),
    \begin{equation*}
      \bm{\delta}_j^{(k)} = \bm{Q}_j - \bm{P}^{(k)}(\bm{t}_j),
      \ j = 0, 1, \cdots, m.
    \end{equation*}
 And then, two procedures are performed, i.e.,
    \emph{vector distribution} and \emph{vector gathering} (Fig.~\ref{fig:lspia_sketch}).
 In the vector distribution procedure,
    all of DVDs corresponding to data points in the $i^{th}$ group
    are distributed to the $i^{th}$ control point $\bm{P}^{(k)}_i$;
 in the vector gathering procedure,
    all of DVDs distributed to the $i^{th}$ control point are weighted averaged to generate the \emph{difference vector for control point} (DVC) (Fig.~\ref{fig:lspia_sketch}),
    \begin{equation*}
      \bm{\Delta}_i^{(k)} =
                \frac{\sum_{j \in I_i} B_i(\bm{t_j}) \bm{\delta}_j}
                     {\sum_{j \in I_i} B_i(\bm{t_j})},
                                 \ i=0, 1, \cdots, n,
    \end{equation*}
    where $I_i$ is the index set of the data points in the $i^{th}$ group.
 Then, the new control point $\bm{P}^{(k+1)}_i$ is produced by adding the
    $i^{th}$ DVC $\bm{\Delta}^{(k)}_i$ to $\bm{P}^{(k)}_i$, i.e.,
    \begin{equation} \label{eq:pkp1}
      \bm{P}^{(k+1)}_i = \bm{P}^{(k)}_i + \bm{\Delta}^{(k)}_i,
      \ i = 0, 1, \cdots, n,
    \end{equation}
    leading to the $(k+1)^{st}$ iteration form,
    \begin{equation} \label{eq:k_to_kp1}
      \bm{P}^{(k+1)}(\bm{t}) = \sum_{i=0}^n \bm{P}_i^{(k+1)} B_i(\bm{t}).
    \end{equation}

 In this way, we get a sequence of iterative forms
    $\{\bm{P}^{k}(\bm{t}),\ k =0,1,\cdots\}$.
 Let,
 \begin{align}
   P^{(k)} & = [\bm{P}^{(k)}_0, \bm{P}^{(k)}_1, \cdots,
                \bm{P}^{(k)}_n]^T, \label{eq:p_matrix}\\
   Q & = [\bm{Q}_0, \bm{Q}_1, \cdots, \bm{Q}_m]^T. \label{eq:q_matrix}
 \end{align}
 From Eq.~\pref{eq:pkp1}, it follows,
 \begin{equation*}
   \begin{split}
        \bm{P}^{(k+1)}_i
        & = \bm{P}^{(k)}_i +
                \frac{1}{\sum_{j \in I_i} B_i(\bm{t}_j)}
                \sum_{j \in I_i} B_i(\bm{t}_j)
                (\bm{Q}_j - \bm{P}^{(k)}(\bm{t}_j)) \\
        & = \bm{P}^{(k)}_i +
                \frac{1}{\sum_{j \in I_i} B_i(\bm{t}_j)}
                \sum_{j \in I_i} B_i(\bm{t}_j)
                \left(\bm{Q}_j -
                \sum_{l=0}^n \bm{P}^{(k)}_l B_l(\bm{t}_j)\right)\\
   \end{split}
 \end{equation*}
 Therefore, we get the LSPIA iterative format in matrix form,
 \begin{equation} \label{eq:it_format}
   P^{(k+1)} = P^{(k)} + \Lambda A^T (Q - A P^{(k)}),
   \ k=0,1,\cdots
 \end{equation}
 where, $\Lambda = diag \left(\frac{1}{\sum_{j \in I_0}B_0(\bm{t}_j)},
                        \frac{1}{\sum_{j \in I_1}B_1(\bm{t}_j)},
                        \cdots,
                        \frac{1}{\sum_{j \in I_n}B_n(\bm{t}_j)} \right)$
 is a diagonal matrix,
 and,
 \begin{equation*}
   A = \begin{bmatrix}
         B_0(\bm{t}_0) & B_1(\bm{t}_0) & \cdots & B_n(\bm{t}_0) \\
         B_0(\bm{t}_1) & B_1(\bm{t}_1) & \cdots & B_n(\bm{t}_1) \\
         \vdots & \vdots &   & \vdots \\
         B_0(\bm{t}_m) & B_1(\bm{t}_m) & \cdots & B_n(\bm{t}_m) \\
       \end{bmatrix}_{(m+1) \times (n+1)}.
 \end{equation*}

 \begin{rem} \label{rem:it_format}
 The iterative format~\pref{eq:it_format} is slightly different from that
    developed in Ref.~\cite{deng2014progressive},
    where diagonal elements of the diagonal matrix $\Lambda$ are equal to each other.
 Although the difference of their iterative formats is slight,
    the convergence analysis of the iterative format~\pref{eq:it_format} is a bit more difficult~\cite{lin2013efficient}.
 \end{rem}

 \begin{rem} \label{rem:digonal_mtx}
 Because diagonal elements of the diagonal matrix
    $\Lambda$ in the iterative format~\pref{eq:it_format} are all positive,
    the diagonal matrix $\Lambda$ is nonsingular.
 \end{rem}

 To show the convergence of the LSPIA iterative format~\pref{eq:it_format},
    it is rewritten as,
    \begin{equation}\label{eq:rewrite_it_fmt}
        P^{(k+1)} = (I - \Lambda A^T A) P^{(k)} + \Lambda A^T Q.
    \end{equation}
 In Ref.~\cite{deng2014progressive},
    it was shown that, when the iterative matrix is \emph{nonsingular},
    the LSPIA iterative format is convergent.
 In the following, we will show that,
    even the matrix $A$ is not of full rank,
    and then $A^T A$ is \emph{singular},
    the iterative format~\pref{eq:it_format} is still convergent.

 We first show some lemmas.

 \begin{lem} \label{lem:nonnegative_eigenvalue}
    The eigenvalues $\lambda$ of the matrix $\Lambda A^T A$ are all real,
    and satisfy $0 \leq \lambda \leq 1$.
 \end{lem}

 \textbf{Proof:} On one hand, suppose $\lambda$ is an arbitrary eigenvalue of the matrix $\Lambda A^T A$ with eigenvector $\bm{v}$, i.e.,
 \begin{equation} \label{eq:eigenvalue}
   \Lambda A^T A \bm{v} = \lambda \bm{v}.
 \end{equation}
 By multiplying $A$ at both sides of Eq.~\pref{eq:eigenvalue}, we have,
 \begin{equation*}
   A \Lambda A^T (A \bm{v}) = \lambda (A \bm{v}).
 \end{equation*}
 It means that $\lambda$ is also an eigenvalue of the matrix $A \Lambda A^T$
    with eigenvector $A\bm{v}$.
 Moreover, $\forall \bm{x} \in R^{m+1}$, because,
 \begin{equation*}
   \bm{x}^T A \Lambda A^T \bm{x} = \bm{x}^T A \Lambda^{\frac{1}{2}} (\Lambda^{\frac{1}{2}})^T A^T \bm{x}
   =  (\bm{x}^T A \Lambda^{\frac{1}{2}}) ( \bm{x}^T A \Lambda^{\frac{1}{2}})^T
   \geq 0,
 \end{equation*}
 the matrix $A \Lambda A^T$ is a positive semidefinite matrix.
 Eigenvalues of a semidefinite matrix are all nonnegative,
    so $\lambda$ is real, and $\lambda \geq 0$.

 On the other hand, because the B-spline basis functions are nonnegative and
    form a partition of unity, it holds, $\norm{A}_{\infty} = 1$.
 Together with $\norm{\Lambda A^T}_{\infty} = 1$, we have,
  \begin{equation*}
    \norm{\Lambda A^T A}_{\infty}  \leq \norm{\Lambda A^T}_{\infty}
    \norm{A}_{\infty} = 1.
  \end{equation*}
 Therefore, the eigenvalue $\lambda$ of matrix $\Lambda A^T A$ satisfies,
 $$\lambda \leq \norm{\Lambda A^T A}_{\infty} \leq 1.$$

 In conclusion, eigenvalues $\lambda$ of the matrix $\Lambda A^T A$ are
    all real,
    and satisfy $0 \leq \lambda \leq 1$.
 $\Box$

 \vspace{0.3cm}

 Because $A^T A$~\pref{eq:rewrite_it_fmt} is singular,
    $\Lambda A^T A$ is also singular,
    and then $\lambda = 0$ is its eigenvalue.
 The following lemma deals with the relationship between the algebraic
    multiplicity and geometric multiplicity of the zero eigenvalue $\lambda = 0$ of $\Lambda A^T A$.

 \begin{rem} \label{rem:rank}
 In this paper, we assume that the dimension of the zero eigenspace of
    $A^T A$ is $n_0$.
 So, the rank of the $(n+1) \times (n+1)$ matrix $A^T A$ is
    $$rank(A^T A) = n-n_0+1.$$
 Because $\Lambda$ is nonsingular (refer to Remark~\ref{rem:digonal_mtx}),
    we have $$rank(\Lambda A^T A) = n-n_0+1.$$
 \end{rem}

 \begin{lem} \label{lem:multiplicity}
    The algebraic multiplicity of the zero eigenvalue of matrix $\Lambda A^T A$ is equal to its geometric multiplicity.
 \end{lem}

\textbf{ Proof:} The proof consists of three parts.

 (1) \emph{The algebraic multiplicity of the zero eigenvalue of matrix $A^T A$ is equal to its geometric multiplicity.}
 Because $A^T A$ is a positive semidefinite matrix,
    it is a diagonalizable matrix.
 Then, for any eigenvalue of $A^T A$ (including the zero eigenvalue),
    its algebraic multiplicity is equal to its geometric multiplicity.
 In Remark~\ref{rem:rank}, we assume that the dimension of the zero
    eigenspace of $A^T A$,
    i.e., the geometric multiplicity of zero eigenvalue of $A^T A$,
    is $n_0$.
 So, the algebraic multiplicity and geometric multiplicity of zero
    eigenvalue of $A^T A$ are both $n_0$.

 \vspace{0.3cm}

 (2) \emph{The geometric multiplicity of the zero eigenvalue of matrix $\Lambda A^T A$ is equal to that of matrix $A^T A$.}
 Denote the eigenspaces of matrices $\Lambda A^T A$ and $A^T A$ associated
    with the zero eigenvalue $\lambda = 0$ as $K_0(\Lambda A^T A)$ and $K_0(A^T A)$, respectively.
 The geometric multiplicities of the zero eigenvalue of matrices
    $\Lambda A^T A$ and $A^T A$ are dimensions of $K_0(\Lambda A^T A)$ and $K_0(A^T A)$, respectively.

 Note that the matrix $\Lambda$ is nonsingular
    (Remark~\ref{rem:digonal_mtx}).
 On one hand, $\forall \bm{v} \in K_0(\Lambda A^T A), \Lambda A^T A \bm{v} = 0$, leading to $A^T A \bm{v} = \Lambda^{-1} 0 = 0$.
 So, $\bm{v} \in K_0(A^T A)$.
 On the other hand, $\forall \bm{w} \in K_0(A^T A), A^T A \bm{w} = 0$,
    resulting in $\Lambda A^T A \bm{w} = 0$.
 So, $\bm{w} \in K_0(\Lambda A^T A)$.
 In conclusion, $K_0(\Lambda A^T A) = K_0(A^T A)$.
 Therefore, the geometric multiplicity of the zero eigenvalue of matrix $\Lambda A^T A$ is equal to that of matrix $A^T A$.

 \vspace{0.3cm}

 (3) \emph{The algebraic multiplicity of the zero eigenvalue of matrix $\Lambda A^T A$ is equal to that of matrix $A^T A$.}
 Denote $I$ as an $(n+1) \times (n+1)$ identity matrix,
    $$A^T A = (b_{ij})_{(n+1) \times (n+1)},\
    \text{and},\ \Lambda = diag(d_0, d_1, \cdots, d_n),$$
    where $d_i > 0, i=0,1,\cdots,n$.

 The characteristic polynomial of $A^T A$ and $\Lambda A^T A$ can be written
    as~\cite[pp.42]{horn1985matrix},
 \begin{equation} \label{eq:char_poly_ata}
   p_{A^T A}(\lambda) = det(\lambda I - A^T A)
   = \lambda^{n+1} - E_1(A^T A) \lambda^n + E_2(A^T A) \lambda^{n-1} + \cdots + (-1)^{n+1}E_{n+1}(A^T A),
 \end{equation}
 and,
 \begin{equation} \label{eq:char_poly_lata}
   p_{\Lambda A^T A}(\lambda) =  det(\lambda I - \Lambda A^T A)
    =  \lambda^{n+1} - E_1(\Lambda A^T A) \lambda^n + E_2(\Lambda A^T A) \lambda^{n-1} + \cdots
    + (-1)^{n+1}E_{n+1}(\Lambda A^T A),
 \end{equation}
 where $E_k(A^T A), k=1,2,\cdots,n+1$ are the sums of the $k \times k$
    principal minors of $A^T A$,
    and $E_k(\Lambda A^T A), k=1,2,\cdots,n+1$ are the sums of the
    $k \times k$ principal minors of $\Lambda A^T A$. 

 On one hand, because the algebraic multiplicity of zero eigenvalue of
    $A^T A$ is $n_0$ (see Part (1)),
    its characteristic polynomial~\pref{eq:char_poly_ata} can be represented as,
    \begin{equation*}
    p_{A^T A} =  \lambda^{n_0} \left( \lambda^{n-n_0+1} - E_1(A^T A) \lambda^{n-n_0} + \cdots + (-1)^{n-n_0+1} E_{n-n_0+1}(A^T A) \right).
    \end{equation*}
    where $E_{n-n_0+1}(A^T A) \neq 0$.
 Moreover, because $A^T A$ is positive semi-definite,
    all of its principal minors are nonnegative.
 Therefore, we have $E_{n-n_0+1}(A^T A) > 0$.
 Consequently, all of $(n-n_0+1) \times (n-n_0+1)$ principal minors of
    $A^T A$ are nonnegative,
    and there is at least one $(n-n_0+1) \times (n-n_0+1)$ principal minor of
    $A^T A$ is positive.

 On the other hand, because $rank(\Lambda A^T A) = n-n_0+1$
    (Remark~\ref{rem:rank}),
    all of $l \times l$ ($l > n-n_0+1$) principal minors of
    $\Lambda A^T A$ are zero.
 Therefore,
 \begin{equation} \label{eq:el}
   E_l(\Lambda A^T A) = 0,\ l > n-n_0+1.
 \end{equation}
 Denote $M_{A^T A}(i_1,i_2,\cdots,i_k)$ and
    $M_{\Lambda A^T A}(i_1,i_2,\cdots,i_k)$ are the $k \times k$ principal minors of $A^T A$ and $\Lambda A^T A$, respectively.
 Now, consider a $k \times k$ principal minor of $\Lambda A^T A$.
 \begin{equation*}
   M_{\Lambda A^T A}(i_1,i_2,\cdots,i_k)  =
               det \begin{pmatrix}
                    d_{i_1} b_{i_1,i_1} & \cdots & d_{i_1} b_{i_1,i_k} \\
                    \cdots & \cdots & \cdots \\
                    d_{i_k} b_{i_k,i_1} & \cdots & d_{i_k} b_{i_k,i_k} \\
               \end{pmatrix}
                = (\prod_{j=1}^k d_{i_j}) M_{A^T A}(i_1,i_2,\cdots,i_k),
 \end{equation*}
 where $\prod_{j=1}^k d_{i_j} > 0$ (Remark~\ref{rem:digonal_mtx}).
 In other words, the principal minor $M_{\Lambda A^T A}(i_1,i_2,\cdots,i_k)$
    of $\Lambda A^T A$ is the product of a principal minor $M_{A^T A}(i_1,i_2,\cdots,i_k)$ of $A^T A$ and a positive value $\prod_{j=1}^k d_{i_j}$.
 Together with that all of $(n-n_0+1) \times (n-n_0+1)$ principal minors of
    $A^T A$ are nonnegative,
    and there is at least one $(n-n_0+1) \times (n-n_0+1)$ principal minor of
    $A^T A$ is positive,
    the sum of all $(n-n_0+1) \times (n-n_0+1)$ principal minors of
    $\Lambda A^T A$, namely, $E_{n-n_0+1}(\Lambda A^T A)$, is positive.
 That is,
 \begin{equation} \label{eq:en}
   E_{n-n_0+1}(\Lambda A^T A) > 0.
 \end{equation}
 By Eqs.~\pref{eq:el} and~\pref{eq:en},
    the characteristic polynomial of $\Lambda A^T A$~\pref{eq:char_poly_lata} can be transformed as,
 \begin{equation*}
    p_{\Lambda A^T A} =  \lambda^{n_0} \left( \lambda^{n-n_0+1} -
        E_1(\Lambda A^T A) \lambda^{n-n_0} + \cdots + (-1)^{n-n_0+1} E_{n-n_0+1}(\Lambda A^T A) \right),
 \end{equation*}
 where $E_{n-n_0+1}(\Lambda A^T A) > 0$.
 It means that the algebraic multiplicity of zero eigenvalue of
    $\Lambda A^T A$ is $n_0$,
    equal to the algebraic multiplicity of zero eigenvalue of $A^T A$.

 \vspace{0.3cm}

 Combing results of part (1)-(3), we have shown that the algebraic
    multiplicity of the zero eigenvalue of matrix $\Lambda A^T A$ is equal to its geometric multiplicity. $\Box$

 \vspace{0.3cm}

 Denote $J_r(a,b)$ as a $r \times r$ matrix block,
 \begin{equation} \label{eq:matrix_block}
    J_r(a,b) = \left(\begin{array}{ccccc}
    a  & b  &   &   & \\
       & a  & b & \text{\huge{0}}  & \\
       &    & \ddots & \ddots & \\
       &  \text{\huge{0}}   &        &   \ddots     &  b \\
       &     &  &  & a\\
 \end{array}\right)_{r \times r}.
 \end{equation}
 Specifically, $J_r(\lambda,1)$ is a $r \times r$ Jordan block.
 Lemma~\ref{lem:nonnegative_eigenvalue} and~\ref{lem:multiplicity} result in
    Lemma~\ref{lem:Jordan} as follows.

 \begin{lem} \label{lem:Jordan}
 The Jordan canonical form of matrix $\Lambda A^T A$~\pref{eq:rewrite_it_fmt}
    can be written as,
    \begin{equation} \label{eq:jordan}
    \small
        J = \left(\begin{array}{ccccccc}
        J_{n_1}(\lambda_1,1) &   &   &   &   &   &  \\
                      & J_{n_2}(\lambda_2,1) &   & \text{\huge{0}} &  &  & \\
                      &   & \ddots &  &  &  &  \\
                      &   &        & J_{n_k}(\lambda_k,1) &  &  &  \\
                      &   &   &           & 0 &  & \\
                      & \text{\huge{0}}  &        &                    &   & \ddots & \\
                      &   &        &                    &   &        & 0 \\
        \end{array}\right)_{(n+1) \times (n+1)},
    \end{equation}
    where $ 0 <\lambda_i \leq 1, i=1,2,\cdots,k$ are nonzero eigenvalues of $\Lambda A^T A$, which need not be distinct, and $J_{n_i}(\lambda_i,1)$~\pref{eq:matrix_block} is an $n_i \times n_i$ Jordan block, $i = 1,2,\cdots,k$.
 \end{lem}

 \textbf{Proof:}
 Based on Lemma~\ref{lem:nonnegative_eigenvalue},
    eigenvalues $\lambda_i$ of $\Lambda A^T A$ are all real and lie in $[0,1]$,
    so the Jordan canonical form of $\Lambda A^T A$ can be written as,
 \begin{equation*}
    \small
        J = \left(\begin{array}{ccccccc}
        J_{n_1}(\lambda_1,1) &   &   &   &   &   &  \\
                      & \ddots  &  &   & \text{\huge{0}} &  &  \\
                      &   &        & J_{n_k}(\lambda_k,1) &  &  &  \\
                      &   &   &           & J_{m_1}(0,1) &  & \\
                      & \text{\huge{0}}  &        &                    &   & \ddots & \\
                      &   &        &                    &   &        & J_{m_l}(0,1) \\
        \end{array}\right)_{(n+1) \times (n+1)},
    \end{equation*}
 where $ 0 <\lambda_i \leq 1, i=1,2,\cdots,k$ are nonzero eigenvalues of $\Lambda A^T A$, which need not be distinct, and $J_{n_i}(\lambda_i,1)$~\pref{eq:matrix_block} is an $n_i \times n_i$ Jordan block, $i = 1,2,\cdots,k$;
 $J_{m_j}(0,1)$ is an $m_j \times m_j$ Jordan block corresponding to the zero eigenvalue of $\Lambda A^T A$, $j=1,2,\cdots,l$.

 According to the theory on Jordan canonical form
    \cite[p.129]{horn1985matrix},
    the number of Jordan blocks corresponding to an eigenvalue is the geometric multiplicity of the eigenvalue,
    and the sum of orders of all Jordan blocks corresponding to an eigenvalue equals its algebraic multiplicity.
 Based on Lemma~\ref{lem:multiplicity},
    the algebraic multiplicity of the zero eigenvalue of matrix
    $\Lambda A^T A$ is equal to its geometric multiplicity,
    so the Jordan blocks $J_{m_j}(0,1),\ j=1,2,\cdots,l$ corresponding to the zero eigenvalue of $\Lambda A^T A$ are all $1 \times 1$ matrix $(0)_{1 \times 1}$.
 This proves Lemma~\ref{lem:Jordan}.
  $\Box$

 \vspace{0.3cm}

 Denote $(A^T A)^+$ as the Moore-Penrose (M-P) pseudo-inverse of the matrix
    $A^T A$.
 We have the following lemma.

 \begin{lem} \label{lem:mp_inverse}
 There exists an orthogonal matrix $V$, such that,
 \begin{equation} \label{eq:mp_inverse}
    V^T (A^T A)^+ (A^T A) V = diag(\underbrace{1,1,\cdots,1}_{n-n_0+1},
                        \underbrace{0,0,\cdots,0}_{n_0}).
 \end{equation}
 \end{lem}

 \textbf{Proof:}
 Because $rank(A^T A) = n-n_0+1$ (Remark~\ref{rem:rank}),
    and $A^T A$ is a positive semidefinite matrix,
    it has singular value decomposition (SVD),
 \begin{equation} \label{eq:svd_decomposition}
   A^T A = V diag(\delta_1, \delta_2, \cdots, \delta_{n-n_0+1},
                \underbrace{0,\cdots,0}_{n_0}) V^T,
 \end{equation}
 where $V$ is an orthogonal matrix, $\delta_i, i = 1,2,\cdots, n - n_0+1$ are singular values of $A^T A$.
 Then, the M-P pseudo-inverse of $A^T A$ is,
 \begin{equation*}
   (A^T A)^+ = V diag(\frac{1}{\delta_1}, \frac{1}{\delta_2}, \cdots, \frac{1}{\delta_{n-n_0+1}}, \underbrace{0,\cdots,0}_{n_0}) V^T.
 \end{equation*}
 Therefore,
 \begin{equation*}
    (A^T A)^+ (A^T A) = V diag(\underbrace{1,1,\cdots,1}_{n-n_0+1},
                        \underbrace{0,0,\cdots,0}_{n_0}) V^T,
 \end{equation*}
  where $V$ is an orthogonal matrix. $\Box$

 \vspace{0.3cm}

 Based on the Lemmas above, we can show the convergence of the iterative
    format~\pref{eq:rewrite_it_fmt} when $A^T A$ is singular.

 \begin{thm} \label{thm:singular}
 When $A^T A$~\pref{eq:rewrite_it_fmt} is singular,
    the iterative format~\pref{eq:rewrite_it_fmt} is convergent.
 \end{thm}

 \textbf{Proof:} By Lemma~\ref{lem:Jordan}, the Jordan canonical form of matrix $\Lambda A^T A$~\pref{eq:rewrite_it_fmt} is $J$~\pref{eq:jordan}.
 Then, there exists a invertible matrix $W$, such that,
 \begin{equation*}
    \Lambda A^T A = W^{-1} J W.
 \end{equation*}
 Therefore (refer to Eq.~\pref{eq:jordan}),
 \begin{equation*}
 \footnotesize
 \begin{split}
   & I - \Lambda A^T A = \\
   & W^{-1}
   \begin{pmatrix}
   J_{n_1}(1-\lambda_1,-1) &   &   &   &   &   &  \\
                      & J_{n_2}(1-\lambda_2,-1) &   & \text{\Large{0}}  &  &  & \\
                      &   & \ddots &  &  &  &  \\
                      &   &        & J_{n_k}(1-\lambda_k,-1) &  &  &  \\
                      &   &        &                    & 1 &  & \\
                      & \text{\Large{0}}  &        &                    &   & \ddots & \\
                      &   &        &                    &   &        & 1 \\
   \end{pmatrix}
   W,
 \end{split}
 \end{equation*}
 where $0 \leq 1-\lambda_i < 1, i = 1,2,\cdots,k$.
 Then, together with Lemma~\ref{lem:mp_inverse}, it holds,
 \begin{equation} \label{eq:it_mtx_converge}
 \begin{split}
    \lim_{l \rightarrow \infty} (I - \Lambda A^T A)^l
     & = W^{-1} diag(\underbrace{0, \cdots, 0}_{n-n_0+1}, \underbrace{1, \cdots, 1}_{n_0}) W \\
     & = I - W^{-1} diag(\underbrace{1,\cdots,1}_{n-n_0+1}, \underbrace{0,\cdots,0}_{n_0}) W \\
     & = I - W^{-1} V^T (A^T A)^+ (A^T A) V W \\
     & = I - (V W)^{-1} (A^T A)^+ (A^T A) (V W)
 \end{split}
 \end{equation}

 Now, consider the linear system $A^T A X = A^T Q$ (refer to Eq.~\pref{eq:q_matrix}).
 It has solutions if and only if~\cite{james1978generalised},
 \begin{equation} \label{eq:solution_condition}
    (A^T A) (A^T A)^+ (A^T Q) = A^T Q.
 \end{equation}

 Subtracting $(A^T A)^+ A^T Q$ from both sides of the iterative
    format~\pref{eq:rewrite_it_fmt},
    together with Eq.~\pref{eq:solution_condition},
    we have,
 \begin{equation} \label{eq:p_k}
 \begin{split}
     & P^{(k+1)} - (A^TA)^+ A^T Q \\
     & = (I - \Lambda A^T A) P^{(k)} +
                                        \Lambda A^T Q - (A^TA)^+ A^T Q \\
       & = (I-\Lambda A^T A) P^{(k)} + \Lambda (A^T A)(A^T A)^+ (A^T Q)
                    - (A^T A)^+ A^T Q \\
       & = (I-\Lambda A^T A) P^{(k)} - (I-\Lambda A^T A)(A^T A)^+ A^T Q\\
       & = (I-\Lambda A^T A) (P^{(k)} - (A^T A)^+ A^T Q) \\
       & = (I-\Lambda A^T A)^{k+1} (P^{(0)} - (A^T A)^+ A^T Q).
 \end{split}
 \end{equation}
 Owing to Eq.~\pref{eq:it_mtx_converge}, it follows,
 \begin{equation} \label{eq:p_infty}
 \begin{split}
    & P^{(\infty)} - (A^TA)^+ A^T Q = \lim_{k \rightarrow \infty} (P^{(k+1)} - (A^TA)^+ A^T Q) \\
    & = \lim_{k \rightarrow \infty} (I-\Lambda A^T A)^{k+1} (P^{(0)} - (A^T A)^+ A^T Q) \\
    & = (I - (V W)^{-1} (A^T A)^+ (A^T A) (V W))
      (P^{(0)} - (A^T A)^+ A^T Q)
 \end{split}
 \end{equation}
 By simple computation, Eq.~\pref{eq:p_infty} changes to,
 \begin{equation} \label{eq:p_infty_simple}
    P^{(\infty)} =  (VW)^{-1} (A^TA)^+ (A^TA) VW (A^TA)^+ A^T Q +
    (I - (V W)^{-1} (A^T A)^+ (A^T A) (V W)) P^{(0)}.
 \end{equation}
 Therefore, the iterative format~\pref{eq:rewrite_it_fmt} is convergent
    when $A^T A$ is singular. Theorem~\ref{thm:singular} is proved. $\Box$

 \vspace{0.3cm}

 \begin{rem}
 Returning to Eq.~\pref{eq:p_infty_simple}, if $V$ is the inverse matrix
    of $W$, i.e., $VW = I$, it becomes,
 \begin{equation} \label{eq:ls_solution}
 \begin{split}
    P^{(\infty)}  & = (A^TA)^+ (A^TA)(A^TA)^+ A^T Q +
    (I - (A^T A)^+ (A^T A)) P^{(0)} \\
     & =(A^TA)^+ A^T Q + (I - (A^T A)^+ (A^T A)) P^{(0)},
 \end{split}
 \end{equation}
 where, $P^{(0)}$ is an arbitrary initial value.
 Eq.~\pref{eq:ls_solution} is the M-P pseudo-inverse solution of the linear
    system $A^T A X = A^T Q$,
    which is the normal equation of the least-squares fitting to the data points~\pref{eq:data_points}.
 Because $P^{(0)}$ is an arbitrary value,
    there are infinite solutions to the normal equation $A^T A X = A^T Q$.
 Within these solutions, $(A^T A)^+ (A^TQ)$ is the one with minimum Euclidean
    norm~\cite{horn1985matrix}.
 \end{rem}

 Actually, if diagonal elements of matrix
    $\Lambda$~\pref{eq:rewrite_it_fmt} are equal to each other,
    denoting as $\alpha$,
    iterative format~\pref{eq:rewrite_it_fmt} can be written as,
    \begin{equation} \label{eq:equal_weight_format}
      P^{(k+1)} = (I - \alpha A^T A) P^{(k)} + \alpha A^T Q.
    \end{equation}
 In this case, we have the following theorem.

 \begin{thm} \label{thm:equal_weight}
    If $A^T A$ is singular,
    and the spectral radius $\rho(\alpha A^T A) \leq 1$,
        the iterative format~\pref{eq:equal_weight_format} converges to the
        M-P pseudo-inverse solution of the linear system $A^T A X = A^T Q$.
    Moreover, if the initial value $P^{(0)} = 0$,
        the iterative format ~\pref{eq:equal_weight_format} converges to $(A^T A)^+ (A^T Q)$, i.e., the M-P pseudo-inverse solution of the linear system $A^T A X = A^T Q$ with the minimum Euclidean norm.
 \end{thm}

 \textbf{Proof:} Because $A^T A$ is both a normal matrix
    and a positive semidefinite matrix,
    its eigen decomposition is the same as its singular value decomposition~\cite{horn1985matrix},
    with the form presented in Eq.~\pref{eq:svd_decomposition}.
 So, we have,
 \begin{equation*}
    \alpha A^T A = V diag(\alpha \delta_1, \alpha \delta_2, \cdots,
                \alpha \delta_{n-n_0+1},
                \underbrace{0,\cdots,0}_{n_0}) V^T,
 \end{equation*}
 where $V$ is an orthogonal matrix, and $\alpha \delta_i, i = 1,2,\cdots, n-n_0+1$ are both the nonzero eigenvalues and nonzero singular values of $\alpha A^T A$.
 Because $\rho(\alpha A^T A) \leq 1$, it holds
    $$0 <\alpha \delta_i \leq 1, i=1,2,\cdots,n-n_0+1.$$
 Then, based on Lemma~\ref{lem:mp_inverse}, we have,
 \begin{equation*}
 \begin{split}
    \lim_{l \rightarrow \infty} (I - \alpha A^T A)^l & =
    V diag(\underbrace{0, \cdots, 0}_{n-n_0+1}, \underbrace{1,\cdots,1}_{n_0}) V^T \\
    & =
        I - V diag(\underbrace{1, \cdots, 1}_{n-n_0+1}, \underbrace{0,\cdots,0}_{n_0}) V^T \\
    & =
        I - V V^T (A^T A)^+ (A^T A) V V^T \\
    &=
        I - (A^T A)^+ (A^T A).
 \end{split}
 \end{equation*}
 Same as the deduction in the proof of Theorem~\ref{thm:singular}
    (Eqs.~\pref{eq:p_k}~\pref{eq:p_infty}), we have,
 \begin{equation*}
 P^{(k+1)} - (A^T A)^+ A^T Q = (I - \alpha A^T A)^{k+1} (P^{(0)} - (A^T A)^+ A^T Q),
 \end{equation*}
 and,
 \begin{equation*}
 \begin{split}
 P^{(\infty)} - (A^TA)^+ A^T Q & = \lim_{k \rightarrow \infty} (P^{(k+1)} - (A^TA)^+ A^T Q) \\
    & = \lim_{k \rightarrow \infty}(I - \alpha A^T A)^{k+1} (P^{(0)} - (A^T A)^+ A^T Q) \\
    & = (I - (A^T A)^+ (A^T A))(P^{(0)} - (A^T A)^+ A^T Q).
 \end{split}
 \end{equation*}
 Therefore,
 \begin{equation*}
 P^{(\infty)} = (A^T A)^+ A^T Q + (I - (A^T A)^+ (A^T A)) P^{(0)},
 \end{equation*}
 where $P^{(0)}$ is an arbitrary initial value.
 It is the M-P pseudo-inverse solution of the linear system
    $A^T A X = A^T Q$,
    and $(A^T A)^+ A^T Q$ is the M-P pseudo-inverse solution with the minimum Euclidean norm. $\Box$

 \section{Example}
 \label{sec:example}

\begin{figure}[!htb]
  \centering
  \subfigure[]{
    \label{subfig:balljoint_input}
    \includegraphics[width=0.3\textwidth]{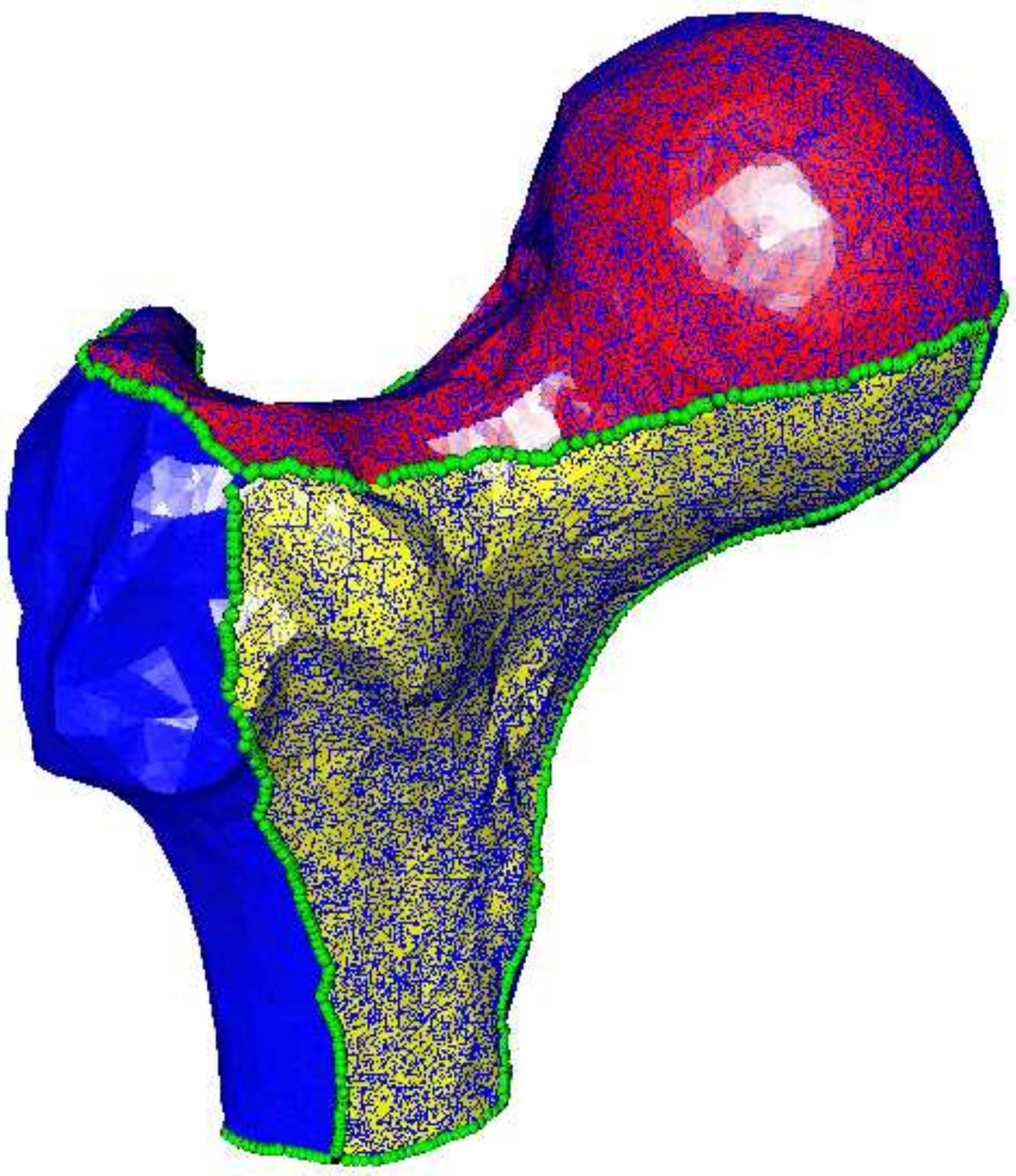}}
  \subfigure[]{
    \label{subfig:balljoint_spline}
    \includegraphics[width=0.3\textwidth]{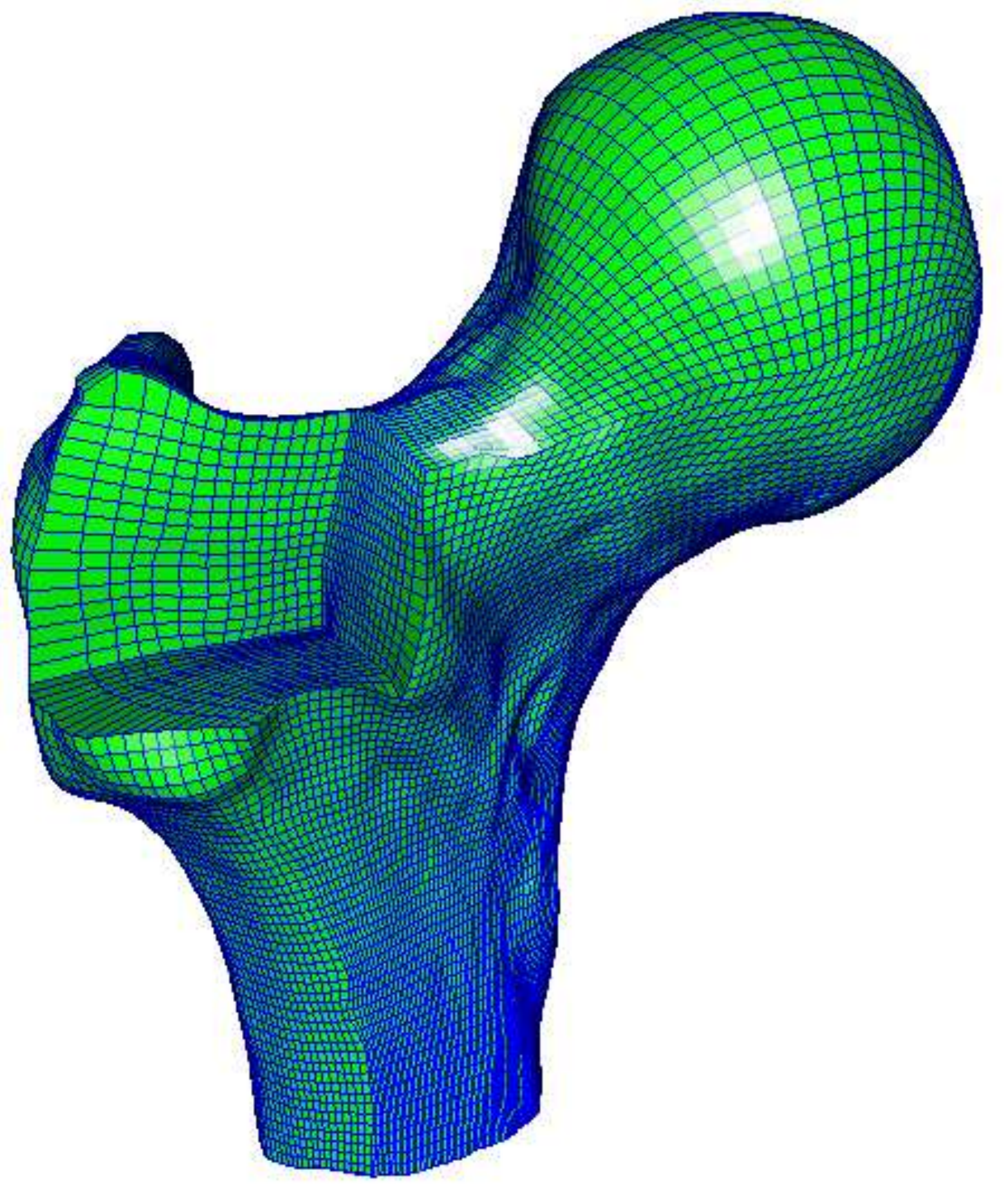}}
    \caption{LSPIA is employed to solve a singular least-squares fitting
        system. The input is a tetrahedral mesh with six segmented
        patches on its boundary (a), and the output is a trivariate B-spline solid (cut-away view) (b).
      }
   \label{fig:example_balljoint}
\end{figure}

 In this section, we present a numerical example,
    which is implemented by Matlab 2013b,
    and run on a PC with $2.9G$ CPU, and $8G$ memory.
 In this example, LSPIA is employed for least-squares fitting a practical
    data point set, i.e.,
    the $43994$ mesh vertices of a tetrahedral mesh model \emph{balljoint} (Fig.~\ref{subfig:balljoint_input}),
     by a tri-cubic trivariate B-spline solid
     (Fig.~\ref{subfig:balljoint_spline}),
    which has $30 \times 21 \times 33$ control points and uniform knot vectors along the three parametric directions with B\'{e}zier end conditions,
    uniformly distributed in the interval $[0,1]$, respectively.
 Fig.~\ref{subfig:balljoint_input} illustrates the input model,
    i.e., a tetrahedral mesh model with six segmented patches on its boundary,
    and Fig.~\ref{subfig:balljoint_spline} is the cut-away view of generated trivariate B-spline solid.
 The tetrahedral mesh vertices are parameterized with the method developed in
    Ref.~\cite{lin2015constructing}.
 In this example, the order of matrix $A^T A$~\pref{eq:rewrite_it_fmt} is
    $20790 \times 20790$,
    and its rank is $12628$.
 Although $A^T A$ is singular, LSPIA converges to a stable solution
    (Fig.~~\ref{subfig:balljoint_spline}) with least-squares fitting error $9.97 \times 10^{-5}$ in $5.50$ seconds.

\section{Conclusions}
\label{sec:conclusion}

 In this paper, we showed that the LSPIA format is convergent when the
    iterative matrix is singular.
 Moreover, when diagonal elements of the diagonal matrix
    $\Lambda$~\pref{eq:it_format} are equal to each other,
    it converges to the M-P pseudo-inverse solution of the least-squares fitting to the given data points.
 Therefore, together with the previous result proved in
    Ref.~\cite{deng2014progressive},
    LSPIA is convergent whatever the iterative matrix is singular or not.
 This property greatly extends the scope of application of LSPIA,
    especially in solving geometric problems in big data processing,
    where singular linear systems frequently appear.

\section*{Acknowledgement} This work is supported by the Natural Science Foundation of China (No. 61379072).



\bibliographystyle{elsart-num}
\bibliography{singular-lsq}

\end{document}